\documentclass{amsart}
\usepackage{amsthm}
\usepackage{amsmath,amssymb}
\usepackage{lipsum}
\usepackage{tikz}
\newtheorem{thm}{Theorem}[section]

\newtheorem{lem}[thm]{Lemma}

\theoremstyle{definition}
\newtheorem{defn}[thm]{Definition}
\theoremstyle{remark}
\newtheorem{rem}[thm]{Remark}
\theoremstyle{example}

\numberwithin{equation}{section}

\newcommand{\R}{\mathbb{R}}

\begin{document}
\title[Lyapunov Type Inequality]
 {Lyapunov type inequality for extremal Pucci's equations}
\author[J.\,Tyagi, R.\,B.\,Verma ]
{J.Tyagi, R.B.\,Verma}
\address{J.\,Tyagi \hfill\break
 Indian Institute of Technology Gandhinagar \newline
 Palaj, Gandhinagar Gujarat, India-382355.}
 \email{jtyagi@iitgn.ac.in, jtyagi1@gmail.com}
\address{R.B.Verma \hfill\break
 Indian Institute of Technology Gandhinagar \newline
 Palaj, Gnadhinagar Gujarat India-382355.}
\email{ram.verma@iitgn.ac.in}
\date{14--06--2017}
\thanks{Submitted 14--06--2017.  Published-----.}
\begin{abstract}
In this article, we establish  Lyapunov type inequality for the following extremal Pucci's equation
\begin{equation*}
\left\{
\begin{aligned}{}
\mathcal{M}^{+}_{\lambda,\Lambda}(D^{2}u)+a(x)u&=0~\text{in}~\Omega,\\
u&=0~\text{on}~\partial\Omega,
\end{aligned}
\right.
\end{equation*}
where $\Omega$ is a smooth bounded domain in $\R^{N},~N\geq2$. This works generalize the well-known works on Lyapunov inequalities to fully nonlinear elliptic equations.
\end{abstract}

\maketitle
\section{Introduction}
The aim of this article is to establish Lyapunov type inequality for the following Pucci's extremal equation:
\begin{equation}\label{coe2}
\left\{
\begin{aligned}{}
\mathcal{M}^{+}_{\lambda,\Lambda}(D^{2}u)+a(x)u&=0~\text{in}~\Omega,\\
u&=0~\text{on}~\partial\Omega,
\end{aligned}
\right.
\end{equation}
where $\Omega$ is a bounded smooth domain in $\R^{N},~N\geq2$. Here, $\mathcal{M}^{+}_{\lambda,\Lambda}$ is called 
Pucci extremal operator. For a given $0<\lambda<\Lambda$, Pucci extremal operator is defined as follows:
\begin{equation}\label{pucci}
\mathcal{M}^{\pm}_{\lambda,\Lambda}(M)=\Lambda\sum_{\pm e_{i}>0}e_{i}+\lambda\sum_{\pm e_{i}<0}e_{i},\\
\end{equation}
where $M$ is an $N\times N$ real symmetric matrix. Let us recall the earlier developments on this subject. The classical Lyapunov 
type inequality says that the necessary condition to have a nontrivial solution to the following boundary value problem
\begin{equation}\label{lya1}
\left\{
\begin{aligned}{}
u''+a(x)u&=0~\text{in}~(c,d),\\
u(c)=u(d)&=0,
\end{aligned}
\right.
\end{equation}
is
\begin{equation}\label{mainf}
\int\limits_{c}\limits^{d}|a|dx>\frac{4}{d-c},
\end{equation}
see \cite{olya, lsurvey}. Later, A. Winter \cite{AW} improved this inequality by replacing $a$ by $a^{+}=\max\{a,0\}$ and also proved that $4$ 
is the optimal constant. Further, this result was generalised to differential equations containing the term $u'$ as well as qusilinear equations, 
see \cite{phertman} and \cite{pode2}, respectively. In \cite{pode2}, Pinasco considered the following one dimensional $p$-Laplace boundary value problem
\begin{equation}\label{coe111}
\begin{aligned}{}
\big(|u'(x)|^{p-2}u'(x)\big)+r(x)|u(x)|^{p-2}u(x)&=0~\text{in}~(c,d),\\
u(c)&=0~\text{and}~u(d)=0,
\end{aligned}
\end{equation}
where $p>1$ and $r$ is a bounded positive function and proved the Lyapunov inequality as well as the lower bound for the eigenvalue as an application of 
Lyapunov inequality. For the recent work in this direction, see \cite{lowerb}. In \cite{pode}, Elbert considered \eqref{coe111} with $r\equiv1$ and proved 
that the necessary condition for \eqref{coe111} with $r\equiv1$ to have a nontrivial solution is
\begin{equation}
\int^{d}_{c}|a|dx~>~\frac{4}{(d-1)^{p-1}}.
\end{equation}
Further, Lee et al.\,\cite{weightp} considered more general operator than a $p$-Laplace operator. In fact, they replace the $p$-Laplace operator in \eqref{coe111} 
with the following operator
\begin{equation}
\big(s(x)|u'(x)|^{p-2}u'(x)\big)+r(x)|u(x)|^{p-2}u(x)=0~\text{in}~(c,d),
\end{equation}
where $s$ is positive, integrable and $r$ is integrable function on $[c,d]$ and established the Lyapunov type inequality. Finally, in this direction, we would like 
to mention the works of  de N\'{a}poli and Pinasco\,\cite{phiode}, where the authors considered the $\psi$-Laplace operator. In fact, they considered the following operator
\begin{equation}
\begin{aligned}{}
(\psi(u'))'+r(x)\psi(u)=0~\text{in}~(c,d),
\end{aligned}
\end{equation}
where $\psi :\R\longrightarrow \R$ is an odd nondecreasing function, such that $\phi(s) = s.\psi(s)$ and $\phi$ is a convex function, see \cite{psi} for similar kind of works. Lyapunov inequality has also been generalised in the context of fractional differential equations, see\cite{frac4,Frac3, Frac4, Frac5, Frac1, Frac2}.\\
We remark that there are interesting works on the Lyapunov inequality for partial differential equation (in short, PDE). 
In fact, the establishment of Lyapunov type inequality for PDE was actively started by generalising the corresponding result for the ODE in \cite{ALP}. 
 Ca\~{n}ada et al.\,\cite{ALP} considered the Neumann problem corresponding to \eqref{lya1}, that is,
\begin{equation}\label{lya22}
\left\{
\begin{aligned}{}
u''+a(x)u&=0~\text{in}~(c,d),\\
u'(c)=u'(d)&=0,
\end{aligned}
\right.
\end{equation}
and defined the following set
\[A=\{a\in L^{1}(c,d)\setminus\{0\}~|~\int^{d}_{c}a(x)dx\geq0~\text{and}~\eqref{lya22}~\text{has non trivial solution}\}.\]
They also defined the following quantity
\begin{equation}\label{beta}
\beta_{p}=\inf_{A\cap L^{p}(c,d)}\|a^{+}\|_{L^{p}(c,d)},
\end{equation}
and studied the qualitative properties of $\beta_{p}$ and obtained the explicit expression for $\beta_{p}$ as a function of $p, b$ and $c.$ There are a good number of applications of Lyapunov inequality, see for instance \cite{sta2,olya2,lsurvey,AJS,SAS,phertman, eib}. The result in \cite{ALP} has been generalised in the context of PDE in \cite{ALF}. More precisely, the authors considered the following problem
 \begin{equation}\label{lya}
\left\{
\begin{aligned}{}
\Delta u+a(x)u&=0~\text{in}~\Omega,\\
u&=0~\text{on}~\partial \Omega,
\end{aligned}
\right.
\end{equation}
where $\Omega\subset\R^{N}(N\geq2)$ is a bounded and smooth domain and the function $a:\Omega\longrightarrow\R,$ belongs to the set $A$ defined as follows:
\[A=\{a\in L^{\frac{N}{2}}(\Omega)\setminus\{0\}~|~\eqref{lya}~\text{has~nontrivial~solutions}\},~\text{if}~N\geq3,~\text{and}\]
\[A=\{\exists~q\in (1,\infty]~\text{with}~q\in L^{q}(\Omega)~\text{and}~\eqref{lya}~\text{has nontrivial solutions}\},~\text{if}~N=2.\]
They also defined the quantity similar to $\beta_{p}$ as in \eqref{beta}, that is,
\[\beta_{p}\equiv\inf_{a\in A\cap L^{p}(\Omega)}\|a^{+}\|_{L^{p}(\Omega)},~\text{for}~1\leq p\leq\infty,\]
and studied the properties of $\beta_{p}.$ More precisely, they proved the following theorem.
\begin{thm}
The following statements hold:\\
\begin{enumerate}
\item[(1)]~If $N=2,$ then $\beta_{p}>0~\iff~1<p\leq\infty.$ If $N\geq3,$ then $\beta_{p}>0~\iff~N/2\leq p\leq\infty$.
\item[(2)]If $N/2<p\leq\infty,$ then $\beta_{p}$ is attained. In this case, any function $a\in A\cap L^{p}(\Omega)$ on which $\beta_{p}$ is attained is of the form:\vspace{0.2cm}
\item[(i)]~$a(x)\equiv\lambda_{1}$, $p\equiv\infty$, where $\lambda_{1}$ is the first eigenvalue.
\item[(ii)]~$a(x)\equiv |u(x)|^{2/(p-1)},$ if $N/2<p<\infty;$ where $u$ is a solution of the problem
\begin{equation*}
\left\{
\begin{aligned}{}
-\Delta u&=|u(x)|^{\frac{2}{p-1}}u(x)~\text{in}~\Omega,\\
u&=0~\text{on}~\partial \Omega.
\end{aligned}
\right.
\end{equation*}
\item[(3)]~The mapping $(N/2,\infty)\longrightarrow\R,$~$p\rightarrow\beta_{p}$ is continuous and the mapping $[N/2, \infty)\longrightarrow\R,$ defined by
\[p\longrightarrow|\Omega|^{-\frac{1}{p}}\beta_{p},\]
is strictly increasing.
\item[(4)] There exists always the limits $\lim_{p\rightarrow\infty}\beta_{p}$ and $\lim_{p\rightarrow(N/2)^{+}}\beta_{p}$ and takes the following values:
\begin{enumerate}
\item[(i)]$\lim_{p\rightarrow\infty}\beta_{p}=\beta_{\infty},$~if $N\geq2;$
\item[(ii)]$\lim_{p\rightarrow(N/2)^{+}}\beta_{p}\geq\beta_{N/2}>0$, if $N\geq3$\\
$\lim_{p\rightarrow1^{+}}\beta_{p}=0,$ if $N=2.$
\end{enumerate}
\end{enumerate}
\end{thm}
In \cite{ALF}, attainability question in the case $\beta_{\frac{N}{2}}$ (i.e, critical case) was left open. This question was settled in \cite{Japa}, by showing that it does not attain. While in contrast to Dirichlet boundary case, in  Neumann boundary case for $N\geq4,$ it is attained, see \cite{minlyap}. Further, the results of \cite{ALF} have been extended to the $p$-Laplace operator with Robin boundary condition in \cite{Robin}. More precisely, they considered the following boundary value problem
\begin{equation}
\left\{
\begin{aligned}{}
-\Delta_{p}u+|u|^{p-2}u&=0~\text{in}~\Omega,\\
|\triangledown u |^{p-2}\frac{\partial u}{\partial\nu}&=a(x)|u|^{p-2}u~\text{on}~\partial\Omega,
\end{aligned}
\right.
\end{equation}
and proved the similar results as in \cite{ALF}. For the Lyapunov type inequality to $p$-Laplace operator with the Dirichlet boundary condition, we refer to \cite{LFA}.\\
Motivated by the above mentioned research and recent works on 
fully nonlinear elliptic equations, see \cite{lu,LMM,AF1,qs,BS,jtv1,jtv2},  there is a natural question to ask.\\
\textbf{Question:}~Can we establish Lyapunov type inequality for fully nonlinear elliptic equations?\\
The aim of this article is to answer this question. 
 More precisely, we establish Lyapunov type inequality for \eqref{coe2}. We remark that the techniques used in earlier research works are not applicable due to the non-divergence 
nature of the problem under consideration. Here, we use another notion of the weak solution,  so-called $L^{N}$-viscosity solution.
We employ Aleksandrov-Bakelman-Pucci estimate for viscosity solutions to get the desired results.
For the definition of $L^{N}$-viscosity solution, 
see Definition \ref{vis}. In order to formulate our results, let us introduce some notations. Let us define $A_{F}$ as follows:
\[A_{F}=\{a\in L^{N}(\Omega)\setminus\{0\}~|~\eqref{coe2}~\text{has~nontrivial~solutions}\}~\text{if}~N\geq2,\]
and set
\begin{equation}
\beta^{F}_{p}=\inf_{a\in A_{F}\cap L^{p}(\Omega)}\|a^{+}\|_{L^{p}(\Omega)}.
\end{equation}
The main result of this paper is the following theorem which we will prove in the next section.
\begin{thm}\label{ply}
The following statements hold:
\begin{enumerate}
\item[(i)]~$\beta^{F}_{p}>0$ for $N\leq p<\infty.$
\item[(ii)]~$\beta^{F}_{p}=0$ for $1\leq p<\frac{N}{2}$, $N\geq 3$ and $\beta^{F}_{1}=0$ for $N=2.$
\item[(iii)] $\beta^{F}_{p}$ is not positive, in general, for $1\leq p<N.$
\end{enumerate}
\end{thm}
We prove Theorem \ref{ply}(iii) through an example. This example also suggests that if we remove a specific class of functions, then we get $\beta^{F}_{p}>0$ for 
all $1\leq p\leq\infty$. For these specific class of functions, see Remark \ref{coer1}. Let us consider the following sets\\
\[P_{g}=\{a\in L^{N}(\Omega)~|~a^{+}(x)>1~\text{for}~a.e~x\in\Omega\},\]
\[P_{l}=\{a\in L^{N}(\Omega)~|~0\leq a^{+}(x)\leq 1~~\text{for}~a.e~x\in\Omega\},\]
and set
\[\tilde{A}_{F}=A_{F}\cap\big(P_{g}\cup P_{l}\big).\]
It is clear that $P_{g}\cap P_{l}=\phi,$ so if $a\in\tilde{A}_{F},$ then either $a\in A_{F}\cap A_{g}$ or $a\in A_{F}\cap A_{l}$.
Now, for $p\geq 1$, set
\[\tilde{\beta}^{F}_{p}=\inf_{a\in\tilde{A}_{F}\cap L^{p}(\Omega)}\|a^{+}\|_{L^{p}(\Omega)},\]
and we also prove the following:
\begin{thm}\label{coep}
For $p\geq 1,$ we have $\tilde{\beta}^{F}_{p}>0$.
\end{thm}
The organisation of this paper is as follows: In Section \ref{basic}, we present important auxiliary results which are used in this article. Section \ref{coe00} is devoted to the proof of main Thoerem \ref{ply} while Section \ref{coep11} contains the proof of Theorem \ref{coep}.
\section{Auxiliary Results and Statements}\label{basic}
We begin this section by recalling the definition of Pucci's extremal operator. For given $0<\lambda<\Lambda$, Pucci extremal operator is defined as follows:
\begin{equation}\label{pucci}
\begin{aligned}{}
\mathcal{M}^{+}_{\lambda,\Lambda}(M)=&\Lambda\sum_{e_{i}>0}e_{i}+\lambda\sum_{e_{i}<0}e_{i},\\
\mathcal{M}^{-}_{\lambda,\Lambda}(M)=&\lambda\sum_{e_{i}>0}e_{i}+\Lambda\sum_{e_{i}<0}e_{i},
\end{aligned}
\end{equation}
where $M$ is a symmetric matrix of size $N\times N$. In general, it is very difficult to find the eigenvalues of the Hessian of a function. But if the given function is radial, that is, there is some $\tilde{u}~:~[0, \infty)\longrightarrow\R$ such that $u(x)=\tilde{u}(|x|)$, then the eigenvalues of the Hessian are given by the following lemma.
\begin{lem}[Lemma 3.1\cite{AF1}]\label{radial}
Let $\tilde{u}~:~[0, \infty)\longrightarrow\R$ be $C^{2}$ function such that $u(x)=\tilde{u}(|x|)$. Then for any $x\in\R^{N}\setminus\{0\}$ the eigenvalues of the Hessian $D^{2}u(x)$ are $\frac{\tilde{u}'(|x|)}{|x|}$ with multiplicity $N-1$ and $\tilde{u}''(|x|)$ with multiplicity 1.
\end{lem}
\begin{defn}\label{vis}
A function $u\in C(\bar{\Omega})\longrightarrow\R,$ is called $L^{N}$-viscosity subsolution (resp. supersolution) to \eqref{coe2} in $\Omega$ if for any $\phi\in W^{2,N}_{loc}(\Omega)$ and any point $x\in\Omega$ at
which $u-\phi$ has local maximum (resp. minimum), we have
\begin{equation}\label{1}
\left\{
\begin{aligned}{}
ess\liminf_{y\rightarrow x}(-\mathcal{M}^{+}_{\lambda,\Lambda}(D^{2}\phi)-a(y)u)&\leq0,\\
(\text{resp.,}~(ess\limsup_{y\rightarrow x}(-\mathcal{M}^{+}_{\lambda,\Lambda}(D^{2}\phi)-a(y)u)&\geq0)).
\end{aligned}
\right.
\end{equation}
\end{defn}
In the proof of our results, Aleksandrov-Bakelman-Pucci (in short, ABP) estimate for viscosity solutions plays an important role. In the context of the viscosity solution this result first of all was proved by Luis A. Caffarelli in the context of continuous viscosity solution, see \cite{lu} and in the context of $L^{N}$-viscosity solution it appears in \cite{LMM}. Further, this result has been generalised in many ways. Here, we adopt ABP estimate from [\cite{BS}, see Theorem 3]. In order to state the theorem, let us set \[\Omega^{\pm}=\{x\in\Omega~:~\pm u(x)>0\}.\]
\begin{thm}\label{abp}
Suppose $u\in C(\bar{\Omega})$ is an $L^{N}$-viscosity solution of
\[\mathcal{M}^{+}_{\lambda,\Lambda}(D^{2}u)\geq f(x)~\text{in}~~(resp.~\mathcal{M}^{-}_{\lambda,\Lambda}(D^{2}u)\leq f(x)),\]
in $\Omega^{+}$ (resp. $\Omega^{-}$) where $f\in L^{N}(\Omega).$ Then
\begin{equation}
\sup_{\Omega}u\leq\sup_{\partial\Omega}u^{+}+\text{diam}(\Omega). C_{1}\|f^{-}\|_{L^{N}(\Omega^{+})}
\end{equation}
(resp. $\sup_{\Omega} u^{-}\leq\sup_{\partial\Omega}u^{-}+\text{diam}(\Omega). C_{1}\|f^{+}\|_{L^{N}(\Omega^{-})}$), where $C_{1}$ is a positive constant which depends on $N,~\lambda, \Lambda, \text{diam}(\Omega)$.
\end{thm} 
\section{Proof of main Theorem}\label{coe00}
\textbf{Proof of Theorem \ref{ply}(i):}
Let us take an arbitrary $a\in A_{F}\cap L^{p}(\Omega)$ and $u$ be a corresponding nontrivial solution to \eqref{coe2}, that is, $u$ satisfies
\begin{equation*}
\begin{aligned}{}
\mathcal{M}^{+}_{\lambda,\Lambda}(D^{2}u)+a(x)u&=0~\text{in}~\Omega,~~\text{i.e,}\\
\mathcal{M}^{+}_{\lambda,\Lambda}(D^{2}u)-a^{-}(x)u&=-a^{+}(x)u~\text{in}~\Omega.\\
\end{aligned}
\end{equation*}
This implies that $u$ satisfies
\[\mathcal{M}^{+}_{\lambda,\Lambda}(D^{2}u)\geq-a^{+}(x)u~\text{in}~\Omega^{+},\]
so by Theorem \ref{abp}, we get
\begin{equation*}
\begin{aligned}{}
\sup_{\Omega}u&\leq\sup_{\partial\Omega}u^{+}+C.~\text{diam}(\Omega)\|a^{+}u\|_{L^{N}(\Omega^{+})},\\
&\leq C.\text{diam}(\Omega)\sup_{\Omega^{+}}u\|a^{+}\|_{L^{N}(\Omega^{+})},~~(\text{since}~u=0~\text{on}~\partial\Omega)\\
&\leq C.\text{diam}(\Omega)\sup_{\Omega}u\|a^{+}\|_{L^{N}(\Omega^{+})},\\
&\leq C.\text{diam}(\Omega)\sup_{\Omega}u\|a^{+}\|_{L^{N}(\Omega)}.\\
\end{aligned}
\end{equation*}
Since $u$ is a nontrivial solution to \eqref{coe2}, so
\[1\leq C.\text{diam}(\Omega)\|a^{+}(x)\|_{L^{N}(\Omega)},~\text{or}\]
\begin{equation}\label{coe3}
\frac{1}{C.\text{diam}(\Omega)}\leq\|a^{+}(x)\|_{L^{N}(\Omega)}.
\end{equation}
Now, if $p=N,$ then by taking the infimum for all $a\in A_{F}\cap L^{N}(\Omega),$ we get the required result, that is,
\[\beta^{F}_{N}\geq\frac{1}{C.\text{diam}(\Omega)}>0.\]
Now on the other hand, if $p>N,$ then
\[\|a^{+}\|_{L^{N}(\Omega)}\leq |\Omega|^{\frac{1}{N}-\frac{1}{p}}\|a^{+}\|_{L^{p}(\Omega)}.\]
So again by \eqref{coe3}, we get
\[\frac{1}{C.\text{diam}(\Omega)}\leq|\Omega|^{\frac{1}{N}-\frac{1}{p}}\|a^{+}\|_{L^{p}(\Omega)},~\text{i.e,}\]
\begin{equation}\label{coe4}
\frac{1}{C.\text{diam}(\Omega)|\Omega|^{\frac{1}{N}-\frac{1}{p}}}\leq\|a^{+}\|_{L^{p}(\Omega)}.
\end{equation}
Since \eqref{coe4} is true for any $a\in A_{F}\cap L^{p}(\Omega),$ so by taking infimum, we again get required result.
\begin{rem}\label{coer2}
In the above proof, we have used the fact that $\Omega^{+}$ is nonempty. However, if $u$ is negative then $\Omega^{+}=\phi$. In this case, we define a function $v=-u$, then $v$ satisfies following equation
\begin{equation}
\left\{
\begin{aligned}{}
\mathcal{M}^{-}_{\lambda,\Lambda}(D^{2}v)+a(x)v&=0~\text{in}~\Omega,\\
v&=0~\text{on}~\partial\Omega,
\end{aligned}
\right.
\end{equation}
for the details, see Remark 2.14 \cite{LMM}. Here $v>0$ is positive so the set $\Omega^{+}_{v}=\{x\in\Omega~|~v(x)>0~\}$ is nonempty, in fact, in this case $\Omega^{+}_{v}=\Omega$. Note also that $v$ satisfies
\begin{equation*}
\begin{aligned}{}
\mathcal{M}^{-}_{\lambda,\Lambda}(D^{2}v)+a(x)v&=0~\text{in}~\Omega,~\text{i.e,}\\
\mathcal{M}^{-}_{\lambda,\Lambda}(D^{2}v)-a^{-}(x)v&=-a^{+}(x)v~\text{in}~\Omega,~~\text{or}\\
\mathcal{M}^{-}_{\lambda,\Lambda}(D^{2}v)&\geq-a^{+}(x)v~\text{in}~\Omega_{v}^{+}.
\end{aligned}
\end{equation*}
Now, using $\mathcal{M}^{-}_{\lambda,\Lambda}(M)\leq \mathcal{M}^{+}_{\lambda,\Lambda}(M)$ for any symmetric matrix $M,$ so we find that $v$ satisfies the following inequality:
\begin{equation}
\begin{aligned}{}
\mathcal{M}^{+}_{\lambda,\Lambda}(D^{2}v)&\geq-a^{+}(x)v~\text{in}~\Omega_{v}^{+},
\end{aligned}
\end{equation}
in $L^{N}$-viscosity sense. Now, repeating the same arguments as in (i), we obtain the required result.
\end{rem}
\textbf{Proof of Theorem \ref{ply} (ii).}
The proof is based on the construction of an example. Here, this example is a modification of an example given for linear case, see[ Lemma 3.1, \cite{ALF}]. First of all, note that if we define $\Omega+x_{0}=\{x+x_{0}~:~x\in\Omega\}$ (for arbitrary $x_{0}\in\R^{N}$), then $\beta^{F}_{p}(\Omega+x_{0})=\beta^{F}_{p}(\Omega)$. On the other hand, if we define $r\Omega=\{rx~:~x\in\Omega\}$ (for arbitrary $r\in\R^{+}$), then $\beta^{F}_{p}(r\Omega)=r^{N/p-2}\beta^{F}_{p}(\Omega)$. Hence
\begin{equation}\label{eee}
\beta^{F}_{p}(\Omega)=0~~\iff~~\beta^{F}_{p}(r\Omega+x_{0})=0.
\end{equation}
Further, for $N\geq 2,$ let us define two numbers
\[\alpha=\frac{\lambda}{\Lambda}(N-1)+1,~\text{and}~\beta=\frac{\Lambda}{\lambda}(N-1)+1,\]
which are frequently used in the construction of examples, below. The proof is divided into two cases separately;  $N\geq3$ and $N=2.$\\
\textbf{Case $N\geq3$.} In view of \eqref{eee}, without loss of generality, we can suppose that $\bar{B}(0,2)\subset\Omega.$ Let us take two arbitrary real numbers $c>d>0$ satisfying $c+d=\frac{\Lambda}{\lambda}(N-1)-1$ and choose $0<\epsilon<\big(\frac{c}{d}\big)^{\frac{1}{d-c}}$. Define the following radial function:
\[
    u(x)=
\begin{cases}
k_{1}|x|^{2}+k_{2},& \text{if }~~|x|\leq \epsilon,\vspace{0.2cm}\\
c|x|^{-d}-d|x|^{-c},& \text{if }~~\epsilon<|x|<1,\vspace{0.2cm}\\
\begin{cases}
\frac{c-d}{(2^{2-\alpha}-1)}[2^{2-\alpha}-|x|^{2-\alpha}],\,\,&\,\text{if}\,\,\alpha<2,\vspace{0.2cm}\\
\frac{d-c}{\log2}[\log\frac{|x|}{2}],& \text{if }~~\alpha=2,\vspace{0.2cm}\\
\frac{c-d}{(1-2^{2-\alpha})}[|x|^{2-\alpha}-2^{2-\alpha}],\,\,&\,\text{if}\,\,\alpha>2
\end{cases}
&\text{if}~~1\leq~|x|\leq 2,\vspace{0.2cm}\\
0,\,\,&\,\,\text{if}\,\,\Omega\cap\{x~|~|x|>2\}.
\end{cases}
\]
In the above expression of $u$, $k_{1}$ and $k_{2}$ are given as follows:
\[k_{1}=\frac{cd}{2}(\epsilon^{-c-2}-\epsilon^{-d-2}),~~\text{and}~~k_{2}=c\epsilon^{-d}\big(1+\frac{d}{2}\big)-d\epsilon^{-c}\big(1+\frac{c}{2}\big).\]
Note that, for $0<\epsilon\leq\Big(\frac{c\big(1+\frac{d}{2}\big)}{d\big(1+\frac{c}{2}\big)}\Big)^{\frac{1}{d-c}},$ $k_{2}\geq0.$ By noting that $c>d,$ it is easy to observe that the following holds:
\begin{equation}\label{obs}
\Big(\frac{c\big(1+\frac{d}{2}\big)}{d\big(1+\frac{c}{2}\big)}\Big)^{\frac{1}{d-c}}>\Big(\frac{c}{d}\Big)^{\frac{1}{d-c}}.
\end{equation}
Now, since $0<\epsilon<\big(\frac{c}{d}\big)^{\frac{1}{d-c}}$ so in view of \eqref{obs}, the functions  $k_{1}$ and $k_{2}$ are positive.
It is easy to see that $u$ is a continuous function and in view of Lemma \ref{radial}, (as in Lemma3.1 \cite{ALF}), it satisfies \eqref{coe2}, where $a$ is given by the following expression
\[ a(x)=
\begin{cases}
\Big(\frac{-2k_{1}N\Lambda}{k_{1}|x|^{2}+k_{2}}\Big),& \text{if }~~|x|\leq \epsilon,\vspace{0.2cm}\\
\frac{\lambda cd}{|x|^{2}},& \text{if }~~\epsilon<|x|<1,\vspace{0.2cm}\\
0,\,\,&\text{if}\,\,x\in\Omega\cap\{x~|~|x|\geq 1\}.
\end{cases}
\]
Here, obviously $a\in L^{\infty}(\Omega)$ and so $a\in A_{F}$ and
\[
    a^{+}(x)=
\begin{cases}
0& \text{if }~~|x|\leq \epsilon,\vspace{0.2cm}\\
\frac{\lambda cd}{|x|^{2}}& \text{if }~~\epsilon<|x|<1,\vspace{0.2cm}\\
0\,\,&\text{if}\,\,x\in\Omega\cap\{x~|~|x|\geq 1\}.
\end{cases}
\]
Now, for $1\leq p<\frac{N}{2}$, let us calculate:
\begin{align}
\|a^{+}\|^{p}_{L^{p}(\Omega)}&=\int\limits_{\Omega}(a^{+})^{p}dx,\nonumber\\
&=\int\limits_{B(0,1)\setminus B(0,\epsilon)}\Big(\frac{\lambda cd}{|x|^{2}}\Big)^{p}dx\nonumber\\
&=\frac{(cd\lambda)^{p}\omega_{N}(1-\epsilon^{N-2p})}{N-2p}.
\end{align}
Thus, we get
\[\|a^{+}\|_{L^{p}(\Omega)}=\Big(\frac{(cd\lambda)^{p}\omega_{N}(1-\epsilon^{N-2p})}{N-2p}\Big)^{\frac{1}{p}}.\]
Now, by definition of $\beta^{F}_{p}$,
\begin{equation}\label{micoe}
\beta^{F}_{p}\leq\Big(\frac{(cd\lambda)^{p}\omega_{N}(1-\epsilon^{N-2p})}{N-2p}\Big)^{\frac{1}{p}}.
\end{equation}
Now, for fixed real numbers $c>d>0$ with $c+d=\frac{\Lambda}{\lambda}(N-1)-1$, we can take limit $\epsilon\rightarrow0,$ in \eqref{micoe} and find that
\begin{equation}\label{micoe1}
\beta^{F}_{p}\leq\frac{cd\omega_{N}^{1/p}}{(N-2p)^{1/p}},
\end{equation}
Finally, taking limit when $d$ tends to zero in \eqref{micoe1}, we conclude that $\beta^{F}_{p}=0,$ and this complete the proof for case $N\geq3$.\vspace{0.2cm}\\
\textbf{Case $N=2$.}~Note that for $N=2,$ for any $x_{0}\in R^{2}$ and $r\in\R^{+},$ we have
  \[\beta^{F}_{p}(r\Omega+x_{0})=\beta^{F}_{p}(\Omega).\]
Thus, again without loss of generality, we can suppose that $\bar{B}(0,2)\subset\Omega.$ Also, for $N=2,$
\[\alpha=\frac{\lambda}{\Lambda}+1~\text{and}~\beta=\frac{\Lambda}{\lambda}+1.\]
Now, take an arbitrary real number $K>\log(\alpha^{2})$ and $\epsilon>0$ satisfying $\log(\epsilon^{2})+K<0$  and consider the following radial function:
\[
    u(x)=
\begin{cases}
\Big(\frac{|x|}{\epsilon}\Big)^{2}+\log\epsilon^{2}+K-1,& \text{if }~~|x|\leq \epsilon,\vspace{0.2cm}\\
\frac{\big[\log(\epsilon^{2}\alpha^{2})|x|^{2-\beta}+\epsilon^{2-\beta}[K-2\log\alpha]-\alpha^{\beta-2}[K+\log\epsilon^{2}]\big]}{[\epsilon^{2-\beta}-\alpha^{\beta-2}]},& \text{if }~~\epsilon<|x|\leq\frac{1}{\alpha},\vspace{0.2cm}\\
-\frac{(1-|x|)^{2}}{\big(1-\frac{1}{\alpha}\big)^{2}}+1+K-\log \alpha^{2},\,\,&\,\text{if}\,\,\frac{1}{\alpha}<|x|\leq1,\vspace{0.2cm}\\
\frac{\big[1+k-\log\alpha^{2}\big]}{\big[2^{2-\alpha}-1\big]}\big[2^{2-\alpha}-|x|^{2-\alpha}\big]\,\,&\text{if}~~1<~|x|\leq 2,\vspace{0.2cm}\\
0,\,\,&\,\,\text{if}\,\,\Omega\cap\{x~|~|x|>2\}.
\end{cases}
\]
As in case $N\geq3$, again it is easy to see that the function $u$ defined above satisfies \eqref{coe2} with $a$ given as follows:
\[
    a(x)=
\begin{cases}
\frac{-4\Lambda}{\big[|x|^{2}+\epsilon^{2}(\log\epsilon^{2}+K-1)\big]},& \text{if }~~|x|\leq \epsilon,\vspace{0.2cm}\\
0,& \text{if }~~\epsilon<|x|<\frac{1}{2},\vspace{0.2cm}\\
\frac{2\Lambda\big[(1\frac{\lambda}{\Lambda})-\frac{1}{|x|}\big]}{\big(1-\frac{1}{\alpha}\big)^{2}\big[-\frac{(1-|x|)^{2}}{(1+\frac{1}{\alpha})^{2}}+K+1-\log\alpha^{2}\big]},\,\,&\,\text{if}\,\,\frac{1}{\alpha}\leq|x|<1,\vspace{0.2cm}\\
0,\,\,&\text{if}~~~|x|\geq 1.
\end{cases}
\]
It is easy to see that $a(x)\geq0$ and $a\in L^{\infty}(\Omega).$ Hence, $a\in A_{F}$. Let us estimate the $L^{1}(\Omega)$ norm of $a:$
\begin{equation}\label{intw}
\begin{aligned}{}
\|a\|_{L^{1}(\Omega)}&=\int_{\Omega}a(x)dx\\
&=\int_{B(0,\epsilon)}a(x)dx+\int_{B(0,1)\setminus B(0,\frac{1}{\alpha})}a(x)dx\\
&=2\pi\int^{\epsilon}\limits_{0}\frac{-4\Lambda r dr}{\big[r^{2}+\epsilon^{2}(\log\epsilon^{2}+K-1)\big]}+\frac{2\pi}{\big(1-\frac{1}{\alpha}\big)}\int^{1}_{\frac{1}{\alpha}}\frac{2\Lambda[\alpha-\frac{1}{r}]rdr}{\big[-\frac{(1-|x|)^{2}}{(1+\frac{1}{\alpha})}+K+1-\log\alpha^{2}\big]}.
\end{aligned}
\end{equation}
The first integral can be evaluated to get the following:
\begin{equation}\label{intw2}
\begin{aligned}
2\pi\int^{\epsilon}\limits_{0}\frac{-4\Lambda r dr}{\big[r^{2}+\epsilon^{2}(\log\epsilon^{2}+K-1)\big]}&=-4\pi\Lambda\big[\log\big(\epsilon^{2}+\epsilon^{2}(\log\epsilon^{2}+k-1)\big)-\log\big(\epsilon^{2}(\log\epsilon^{2}+K-1)\big)\big]\\
&=4\pi\Lambda\log\Big(\frac{\epsilon^{2}\big(\log\epsilon^{2}+K-1\big)}{\epsilon^{2}+\epsilon^{2}\big(\log\epsilon^{2}+K-1\big)}\Big).
\end{aligned}
\end{equation}
The second integral in \eqref{intw} can be estimated as follows:
\begin{equation}\label{inw3}
\begin{aligned}{}
\frac{2\pi}{\big(1-\frac{1}{\alpha}\big)}\int^{1}_{\frac{1}{\alpha}}\frac{2\Lambda[\alpha-\frac{1}{r}]rdr}{\big[-\frac{(1-|x|)^{2}}{(1+\frac{1}{\alpha})}+K+1-\log\alpha^{2}\big]}&=\frac{4\pi\Lambda}{\big(1-\frac{1}{\alpha}\big)}\int^{1}\limits_{\frac{1}{\alpha}}\frac{[\alpha r-1]dr}{\big[-\frac{(1-r)^{2}}{\big(1-\frac{1}{\alpha})^{2}}+K+1-\log\alpha^{2}\big]}\\
&\leq\frac{4\pi\Lambda}{\big(1-\frac{1}{\alpha}\big)}\frac{[\alpha-1]}{\big[-\frac{(1-\frac{1}{\alpha})^{2}}{\big(1-\frac{1}{\alpha})^{2}}+K+1-\log\alpha^{2}\big]}\int^{1}\limits_{\frac{1}{\alpha}}dr\\
&=\frac{4\pi\Lambda}{\big(1-\frac{1}{\alpha}\big)}\frac{[\alpha-1]}{\big[K-\log\alpha^{2}\big]}\big[1-\frac{1}{\alpha}\big]\\
&=\frac{4\pi\Lambda\alpha}{[K-\log\alpha^{2}]}.
\end{aligned}
\end{equation}
On combining \eqref{intw},\eqref{intw2} and \eqref{inw3}, we find that
\[\|a\|_{L^{1}(\Omega)}\leq4\pi\Lambda\log\Big(\frac{\epsilon^{2}\big(\log\epsilon^{2}+K-1\big)}{\epsilon^{2}+\epsilon^{2}\big(\log\epsilon^{2}+K-1\big)}\Big)+\frac{4\pi\Lambda\alpha}{[K-\log\alpha^{2}]}.\]
Thus,
\begin{equation}\label{ra}
\beta^{F}_{1}\leq4\pi\Lambda\log\Big(\frac{\epsilon^{2}\big(\log\epsilon^{2}+K-1\big)}{\epsilon^{2}+\epsilon^{2}\big(\log\epsilon^{2}+K-1\big)}\Big)+\frac{4\pi\Lambda\alpha}{[K-\log\alpha^{2}]}.
\end{equation}
But, for fixed real number $K>\log(\alpha^{2}),$ we can take limit $\epsilon$ tending to zero in \eqref{ra} to get
\begin{equation}\label{intw4}
\beta^{F}_{1}\leq\frac{4\pi\Lambda\alpha}{[K-\log\alpha^{2}]}.
\end{equation}
Finally, taking limit as $K$ approaching to $+\infty,$ we conclude that $\beta^{F}_{1}=0.$\vspace{0.3cm}\\
{\textbf{Proof of Theorem \ref{ply}(iii).}~
We prove this part by constructing a simple example. Let us consider the following problem
\begin{equation}\label{coe112}
\left\{
\begin{aligned}{}
\mathcal{M}^{+}_{\lambda,\Lambda}(D^{2}u)+a(x)u&=0~\text{in}~B(0,2\lambda(N-1)),\\
u&=0~\text{on}~\partial B(0, 2\lambda(N-1)).
\end{aligned}
\right.
\end{equation}
Next onwards, we denote $\lambda(N-1)$ by $\bar{r}$, so $2\lambda(N-1)=2\bar{r}$. Let us also define a number
$\alpha=\frac{\lambda}{\Lambda}(N-1)+1$, and consider the following function:
\[u(x)=
\begin{cases}
    e^{-\frac{\bar{r}}{k(k+1)}},& \text{if }~~|x|<\frac{\bar{r}}{(k+1)},\vspace{0.2cm}\\
    e^{-\frac{|x|}{k}},& \text{if }~~\frac{\bar{r}}{(k+1)}\leq |x|<\frac{3\bar{r}}{2},\vspace{0.2cm}\\
\begin{cases}
\frac{e^{-\frac{3\bar{r}}{2k}}}{[2^{2-\alpha}-(\frac{3}{2})^{2-\alpha}]}\big[2^{2-\alpha}-(\frac{|x|}{\bar{r}})^{2-\alpha}\big],& \text{if } \alpha\not=2\vspace{0.2cm}\\
\frac{e^{-\frac{3\bar{r}}{2k}}}{\log\frac{4}{3}}\big[\log(2\bar{r})-\log |x|\big],& \text{if }~~\alpha=2.\\
\end{cases}
&;\text{if}~~\frac{3\bar{r}}{2}\leq |x|\leq 2\bar{r},\\
\end{cases}
\]
It is easy to verify that for each $k,$ $u$ satisfies \eqref{coe112}, where $a$ is given by
\[a(x)=
\begin{cases}
    0,& \text{if }~~|x|<\frac{\bar{r}}{(k+1)},\vspace{0.2cm}\\
\frac{\bar{r}}{k|x|}-\frac{\Lambda}{k^{2}},& \text{if }~~\frac{\bar{r}}{(k+1)}\leq |x|<\frac{3\bar{r}}{2},\vspace{0.2cm}\\
0,~&\text{if}~~\frac{3\bar{r}}{2}\leq |x|\leq 2\bar{r}.
\end{cases}
\]
Thus, it is clear that for each $k$
\[a^{+}(x)\leq
\begin{cases}
0,& \text{if }~~|x|<\frac{\bar{r}}{(k+1)},\vspace{0.2cm}\\
\frac{\bar{r}}{k|x|},& \text{if }~~\frac{\bar{r}}{(k+1)}\leq |x|<\frac{3\bar{r}}{2},\vspace{0.2cm}\\
0,~&\text{if}~~\frac{3\bar{r}}{2}\leq |x|\leq 2\bar{r}.
\end{cases}
\]
Now, let us assume that $1\leq p<N$, and compute:
\begin{equation*}
\begin{aligned}{}
\|a^{+}\|^{p}_{L^{p}(B(0,2\bar{r}))}&=\int\limits_{B(0,2\bar{r})}|a^{+}(x)|^{p}dx\\
&=\int\limits_{\{\frac{\bar{r}}{k+1}\leq |x|<\frac{3\bar{r}}{2}\}}(a^{+}(x))^{p}dx\\
&\leq\int\limits_{\{\frac{\bar{r}}{k+1}\leq |x|<\frac{3\bar{r}}{2}\}}\Big(\frac{\bar{r}}{k|x|}\Big)^{p}dx\\
&=\omega_{N}\int\limits^{3\bar{r}/2}\limits_{\bar{r}/k+1}r^{N-1}\Big(\frac{\bar{r}}{kr}\Big)^{p}dr\\
&=\omega_{N}\Big(\frac{\bar{r}}{k}\Big)^{p}\int\limits^{3\bar{r}/2}\limits_{\bar{r}/k+1}r^{N-p-1}dr\\
&=\frac{\omega_{N}}{(N-p)}\Big(\frac{\bar{r}}{k}\Big)^{p}\Big[r^{N-p}\Big]^{3\bar{r}/2}_{\bar{r}/k+1}\\
&=\frac{\omega_{N}}{(N-p)}\Big(\frac{\bar{r}^{N}}{k^{p}}\Big)\Big[\Big(\frac{3}{2}\Big)^{N-p}-\frac{1}{(k+1)^{N-p}}\Big].
\end{aligned}
\end{equation*}
Thus, we find that $\|a^{+}\|_{L^{p}(B(0,2\bar{r}))}\longrightarrow0$ as $k\rightarrow\infty.$ Consequently, for $\Omega=B(0,2\lambda(N-1)),$
\[\beta^{F}_{p}=\inf_{a\in A_{F}\cap L^{p}(\Omega)}\|a^{+}\|_{L^{p}(\Omega)}=0.\]
\begin{rem}\label{coer1}
In the above example, we have shown that for each $k,$ the following problem:
\begin{equation}
\left\{
\begin{aligned}{}
\mathcal{M}^{+}_{\lambda,\Lambda}(D^{2}u_{k})+a_{k}(x)u_{k}&=0~\text{in}~B(0,2\lambda(N-1)),\\
u_{k}&=0~\text{on}~\partial B(0, 2\lambda(N-1)),
\end{aligned}
\right.
\end{equation}
has a nontrivial solution and for $1\leq p<N$, $\|a_{k}^{+}\|_{L^{p}(\Omega)}\rightarrow0$ as $k\rightarrow\infty$.\\
Now, choose $\tilde{k}$ large enough such that for $k\geq\tilde{k},$ the following hold:
\begin{enumerate}
\item[(i)]~$\frac{\Lambda}{k}<1,$
\item[(ii)]~$\frac{k}{\Lambda}\geq\frac{3}{2}.$
\end{enumerate}
Note that, if $k\geq\tilde{k},$ then by (ii), we have
\begin{align}\nonumber
\frac{3}{2}&\leq\frac{k}{\Lambda},~\text{i.e},\\ \label{lref}
\frac{3k\bar{r}}{2}&\leq\frac{k^{2}\bar{r}}{\Lambda}.
\end{align}
In view of \eqref{lref}, for any $\frac{\bar{r}}{(k+1)}\leq |x|\leq\frac{3\bar{r}}{2},$ we have
\begin{align}\nonumber
k|x|&\leq\frac{k^{2}\bar{r}}{\Lambda},~~\text{or}\\ \nonumber
\frac{k|x|}{\bar{r}}&\leq\frac{k^{2}}{\Lambda},~~\text{i.e},\\ \label{llref}
\frac{\Lambda}{k^{2}}&\leq\frac{\bar{r}}{k|x|}.
\end{align}
Thus, for any $k\geq\tilde{k},$ definition of $a_{k}$ and \eqref{llref}, yields the following
\[a^{+}_{k}(x)=a_{k}(x).\]
Now for $k\geq\tilde{k}$, consider the following set
\begin{equation*}
\begin{aligned}{}
G_{1}&=\{x\in B(0,2\bar{r})~|~a^{+}_{k}(x)>1\},\\
&=\left\{\frac{\bar{r}}{(k+1)}\leq |x|<\frac{3\bar{r}}{2}~|~\frac{\bar{r}}{k|x|}-\frac{\Lambda}{k^{2}}>1 \right\},\\
&=\left\{\frac{\bar{r}}{(k+1)}\leq |x|<\frac{3\bar{r}}{2}~|~\frac{\bar{r}}{k|x|}>1+\frac{\Lambda}{k^{2}}\right\},\\
&=\left\{\frac{\bar{r}}{(k+1)}\leq |x|<\frac{3\bar{r}}{2}~|~\frac{\bar{r}}{k\big(1+\frac{\Lambda}{k^{2}}\big)}>|x|\right\},\\
&=\left\{\frac{\bar{r}}{(k+1)}\leq |x|<\frac{3\bar{r}}{2}~|~\frac{\bar{r}}{\big(k+\frac{\Lambda}{k}\big)}>|x|\right\}.\\
\end{aligned}
\end{equation*}
Note that, in view of assumption (i), for $k\geq\tilde{k}$, we have $\frac{\Lambda}{k}<1,$ so
\begin{equation}\label{ee11}
\frac{\bar{r}}{\big(k+\frac{\Lambda}{k}\big)}>\frac{\bar{r}}{\big(k+1\big)},~\text{for}~k\geq\tilde{k}.
\end{equation}
Also notice that for any $k\geq1,$ we have
\begin{equation}\label{ee12}
\frac{\bar{r}}{\big(k+\frac{\Lambda}{k}\big)}<\frac{\bar{r}}{k}\leq\bar{r}<\frac{3\bar{r}}{2},~\text{for}~k\geq1.
\end{equation}
Thus, in view of \eqref{ee11} and \eqref{ee12}, $G_{1}$ takes the following form
\[G_{1}=\{x\in\R^{N}~|~\frac{\bar{r}}{(k+1)}\leq |x|<\frac{\bar{r}}{\big(k+\frac{\Lambda}{k}\big)}\}.\]
Hence the Lebesgue measure of $G_{1},$ i.e,
\[|G_{1}|=\bar{r}^{N}\omega(N)\Big[\frac{1}{\big(k+\frac{\Lambda}{k}\big)^{N}}-\frac{1}{\big(k+1\big)^{N}}\Big]\not=0,\]
for $k\geq\tilde{k}.$ Of course the Lebesgue measure of the following set
\[L_{1}=\{x\in B(0, 2\bar{r})~|~0\leq a^{+}_{k}(x)\leq 1\},\]
is not zero.
\end{rem}
In view of the above, it is natural ask that if we remove those functions from $A_{F},$ for which Lebesgue measure $|G_{1}|\not=0$ and $|L_{1}|\not=0$, then, 
whether the modified quantity corresponding to $\beta^{F}_{p}$ is positive or not? In fact, the answer to this question is affirmative. Next, in the proof of Theorem \ref{coep}, 
we answer to this question.
\section{Proof of Theorem \ref{coep}}\label{coep11}
If $N\leq p,$ then from the definition of $\beta^{F}_{p}$, $\beta^{F}_{p}\leq\tilde{\beta}^{F}_{p},$ so in this case, the result follows form Theorem \ref{ply}. Now, we consider 
the case $1\leq p<n$. Let us take an arbitrary $a\in \tilde{A}_{F}$ and let $u$ be a corresponding nontrivial solution to \eqref{coe2}. Now, if $a\in P_{g},$ then for 
$x\in \Omega$ a.e, we have
\begin{equation*}
\begin{aligned}{}
1&<a^{+}(x),\\
1&<(a^{+}(x))^{p}=|a^{+}(x)|^{p},~\text{for any}~1\leq p<N.
\end{aligned}
\end{equation*}
So an integration yields that
\begin{equation}\label{coe5}
\begin{aligned}{}
|\Omega|&<\int\limits_{\Omega}|a^{+}(x)|^{p}dx,~\text{i.e,}\\
|\Omega|^{\frac{1}{p}}&<\|a^{+}\|_{L^{p}(\Omega)}.
\end{aligned}
\end{equation}
Now, otherwise, if $a\in P_{l},$ so we have $0\leq a^{+}(x)\leq 1$ a.e $x\in\Omega$. Since there is a nontrivial solution to
\begin{equation*}
\begin{aligned}{}
\mathcal{M}^{+}_{\lambda,\Lambda}(D^{2}u)+a(x)u&=0~\text{in}~\Omega,\\
\mathcal{M}^{+}_{\lambda,\Lambda}(D^{2}u)-a^{-}(x)u&=-a^{+}(x)u~\text{in}~\Omega.\\
\end{aligned}
\end{equation*}
That is,
\begin{equation}\label{coe6}
\mathcal{M}^{+}_{\lambda,\Lambda}(D^{2}u)\geq-a^{+}(x)u~\text{in}~\Omega^{+}.
\end{equation}
In order to get the required result, we need to adjust the right hand side of \eqref{coe6} before applying the (ABP) estimate. So let us proceed. Since $0\leq a^{+}\leq1,$ so
\[0\leq a^{+}(x)\leq (a^{+}(x))^{q}~\text{for any}~0\leq q<1.\]
In particular, since $\frac{p}{N}<1$ so $0\leq a^{+}(x)\leq (a^{+}(x))^{\frac{p}{N}}$. Now by \eqref{coe6}, we find that
\begin{equation}
\mathcal{M}^{+}_{\lambda,\Lambda}(D^{2}u)\geq-(a^{+}(x))^{\frac{p}{N}}u~\text{in}~\Omega^{+}.
\end{equation}
Noting that $a^{+}\in L^{p}(\Omega)$, and $u\in C(\bar{\Omega}),$ we conclude that $|a^{+}(x)|^{\frac{p}{N}}u\in L^{N}(\Omega)$. Therefore by \eqref{abp}, we get
\begin{align}
\sup_{\Omega}u&\leq\sup_{\partial\Omega}u^{+}+C.\text{diam}(\Omega)\|(a^{+})^{\frac{p}{N}}(x)u\|_{L^{N}(\Omega^{+})} \nonumber\\ \nonumber
&\leq C.\text{diam}(\Omega)\sup_{\Omega^{+}}u\|(a^{+}(x))^{\frac{p}{N}}\|_{L^{N}(\Omega^{+})}~~(\text{since}~u=0~\text{on}~\partial\Omega)\\ \nonumber
&\leq C.\text{diam}(\Omega)\sup_{\Omega}u\|(a^{+}(x))^{\frac{p}{N}}\|_{L^{N}(\Omega^{+})}\\
&\leq C.\text{diam}(\Omega)\sup_{\Omega}u\|(a^{+}(x))^{\frac{p}{N}}\|_{L^{N}(\Omega)}.\label{ran1}
\end{align}
Now
\begin{equation}\label{ran2}
\begin{aligned}{}
\|(a^{+}(x))^{\frac{p}{N}}\|^{N}_{L^{N}(\Omega)}&=\int\limits_{\Omega}|(a^{+}(x))^{\frac{p}{N}}|^{N}
&=\int\limits_{\Omega}|(a^{+}(x))|^{p}
&=\|a^{+}\|^{p}_{L^{p}(\Omega)}.
\end{aligned}
\end{equation}
Now \eqref{ran1} and \eqref{ran2} yield that
\[\sup_{\Omega}u\leq C.\text{diam}(\Omega)\sup_{\Omega}u\|a^{+}\|^{\frac{p}{N}}_{L^{p}(\Omega)},\]
i.e,
\begin{equation}\label{coe7}
\Big(\frac{1}{C.\text{diam}(\Omega)}\Big)^{\frac{N}{p}}\leq\|a^{+}\|_{L^{p}(\Omega)}.
\end{equation}
Now on combining \eqref{coe7} and \eqref{coe5}, we find that
\[\min\Big\{|\Omega|^{\frac{1}{p}}, \Big(\frac{1}{C.\text{diam}(\Omega)}\Big)^{\frac{N}{p}}\Big\}\leq\|a^{+}\|_{L^{p}(\Omega)}.\]
Since $a\in\tilde{A}_{F}\cap L^{p}(\Omega)$ is arbitrary so by taking infimum over $a,$ we find $\tilde{\beta}^{F}_{p}>0$, and this completes the proof.

\end{document}